\documentclass[]{article}
\usepackage{enumerate,color}
\usepackage{amsmath,amsthm}
\usepackage{amsfonts}
\usepackage{amssymb}

\begin{document}
\title{On recurrence and transience of multivariate near-critical stochastic processes}
\author{G\"otz Kersting\thanks{Institut f\"ur Mathematik, Goethe Universit\"at, Frankfurt am Main, Germany, kersting@math.uni-frankfurt.de, work partially supported by the DFG Priority Programme SPP 1590 ``Probabilistic Structures in Evolution''}}
\date{\today}
\maketitle

\begin{abstract} 
We obtain complementary recurrence and transience criteria for processes $X=(X_n)_{n \ge 0}$ with values in $\mathbb R^d_+$ fulfilling
a non-linear equation $X_{n+1}=MX_n+g(X_n)+ \xi_{n+1}$. Here $M$ denotes a primitive matrix  having Perron-Frobenius eigenvalue 1, and $g$ denotes some function. The conditional expectation and variance of the noise $(\xi_{n+1})_{n \ge 0}$ are such that $X$ obeys a weak form of the Markov property. The results generalize criteria for the \mbox{1-dimensional} case in \cite{ke1}. 
\medskip



\end{abstract}

\section{Introduction and main results}

For Markov chains with   a higher-dimensional state space it is  in general difficult to obtain criteria for recurrence or transience which cover a broader class of models. Typically this requires  some specific assumptions on the typ of model. In this paper we consider discrete time stochastic processes $X=(X_n)_{n \ge 0}$ taking values in the positive orthant $\mathbb R^d_+$  (consisting of column vectors) with $d \ge 1$,  which obey non-linear equations of the form
\begin{align}
X_{n+1}=MX_n + g(X_n)+ \xi_{n+1} \ , \quad n \in \mathbb N_0 \ . 
\label{gleichung}
\end{align}
Here $M$ denotes a $d\times d$ matrix with non-negative entries and $g:\mathbb R^d_+ \to \mathbb R^d_+$ a measurable function.  Let us successively discuss our assumptions on $M$, $g$ and the random fluctuations $(\xi_{n+1})_{n\ge 0}$.

We require that $M$ is a primitive matrix meaning that for a certain power of $M$ all entries are (strictly) positive. Then it is known from Perron-Frobenius theory that $M$ has left and right eigenvectors $\ell=(\ell_1, \ldots, \ell_d)$ and $r=(r_1, \ldots, r_d)^T$ belonging to some positive eigenvalue and possessing only  positive entries. We assume that this eigenvalue is 1:
\[ \ell M= \ell \ , \quad Mr=r \ . \]
Further $\ell$ and $r$ are unique up to scaling factors. As is customary we choose them such that
\begin{align} \ell r=1\ . \label{normierung} \end{align}

For the function $g$ we assume that  
\begin{align} \| g(x)\| = o(\| x \|) \text{ as } \| x\| \to \infty  
\label{smallo}
\end{align}
with some norm $\| \ \|$ on the Euclidian space $\mathbb R^d$.

 As to the random fluctuations we demand that $X$ is adapted to a filtration $\mathcal F=(\mathcal F_n)_{n \ge 0}$ such that
\begin{align} \mathbf E[\ell \xi_{n+1} \mid \mathcal F_n]= 0 \ , \quad \mathbf E[(\ell\xi_{n+1})^2 \mid \mathcal F_n] = \sigma^2(X_n) \text{ a.s.} 
\label{variance}
\end{align}
for some measurable function $\sigma: \mathbb R^d_+ \to \mathbb R_+$ fulfilling
\begin{align}
\sigma(x)= o(\| x \|) \quad \text{for } \| x\| \to \infty  \ .
\label{smallo2}
\end{align} 

In view of applications such as branching processes we might summerize these requirements on the whole as the assumption of {\em near criticality}. Quite a few models fit into this framework.  Here we do not dwell on them  but refer to the paper \cite{ke2} and to the literature cited therein. The assumption \eqref{variance} establishes a weak form of the Markov property.  We do not assume that $X$ is a Markov chain but just formulate those assumptions which are required for the martingale considerations in our proofs. Certainly applications of our results will typically concern Markov chains.

The aim of this paper is to establish criteria which allow to decide whether $\|X_n\| \to \infty $ is an event of zero probability or not. Loosely speaking these are criteria for {\em recurrence} or {\em transience} of our models. In the univariate case $d=1$ this question has been discussed in \cite{ke1}. Ignoring some side conditions the result there was as follows:  If for some $\varepsilon >0$ and for $x$ sufficiently large
\[xg(x) \le\frac {1-\varepsilon} 2\sigma^2(x)  \ , \]
then we have recurrence. If on the other hand for some $\varepsilon >0$ and for $x$ sufficiently large
\[ xg(x) \ge \frac {1+\varepsilon}2 \sigma^2(x) \ , \]
then there is transience. Heuristically this can be understood as follows: In the first regime it is the noise $\xi_{n+1}$ which dominates the drift $g(X_n)$, while in the second regime it is the other way round.  We like to generalize this dichotomy to the multivariate setting. 

A possible way of generalization is to suitably convert each of the two conditions to {\em all} $x \in \mathbb R^d_+$ with  sufficiently large norm $\|x\|$, see Klebaner \cite{kleb} and Gonz\'alez et al \cite{go}. A relaxation of this approach for special choices of $g$ and $\sigma^2$ covering new examples has been obtained by Adam \cite{ad}. Yet one can do with weaker assumptions. The intuition behind this assertion is that our processes behave in a sense 1-dimensional. More precisely, if the event $\|X_n\|\to \infty$ occurs, then in view of \eqref{smallo} and \eqref{smallo2} it is the term $MX_n$, which dominates on the right-hand side  of \eqref{gleichung}. Thus one would expect that   $X_n$ will escape to $\infty$  approximately along  the ray $\overline r=\{ \nu r: \nu \ge 0\}$ spanned by the eigenvector $r$ of $M$. This suggests that the two conditions above are required  only in  certain vicinities of this ray. (The last assertion of Theorem 2 below confirms this heuristics.) 

To formalize these considerations let us introduce some notation.
For any $x \in \mathbb R^d$ let
\[ \hat x:= r\ell x \ , \quad \check x:=(I -r\ell)x \ , \quad \text{thus } x=\hat x+ \check x\ , \]
with the identity matrix $I$. Note that $\hat x$ is the multiple $(\ell x)r$ of the vector $r$ and thus belongs to the ray $\overline r$. From \eqref{normierung} $r\ell r\ell=r\ell$ respectively $\hat {\hat x}=\hat x$ meaning that $r\ell$ is a projection matrix.  Moreover  $\ell \hat x=\ell x$ or $\ell \check x=0$. The two conditions 
$\hat x \in \overline r$ and $\ell \hat x=\ell x $ 
determine $\hat x\in \mathbb R^d$  uniquely.

For convenience we reqire the additional moment condition (which could be relaxed)
\begin{align} \exists \delta>0, c <\infty \ : \  \mathbf E[\| \xi_{n+1} \|^p \mid \mathcal F_n] \le c \sigma^p(X_n) \quad \text{with } p=2+\delta\ .
\tag{A1}
\label{A1}
\end{align}
\paragraph{Theorem 1.} {\em Let }(A1) {\em be fulfilled and let $\varepsilon >0$. Assume that  for every $b>0$ there exists some $a>0$ such that for $x \in \mathbb R_+^d$
\begin{align} \| x\| \ge a \ , \ \| \check x\|^2 \le b \|x\| \cdot \| g(x)\| \quad \Rightarrow \quad \ell x \cdot \ell g(x) \le \frac{1- \varepsilon}2 \sigma^2(x) \ . 
\label{driftcond}
\end{align}
Then }
\[ \mathbf P(\| X_n\| \to \infty ) = 0 \ .  \]

\mbox{}\\
In the case $d=1$ we have $\check x=0$ and $\ell x \cdot \ell g(x)= xg(x)$ such that we are back to the result from \cite{ke1}. Note that due to \eqref{smallo} the above condition $ \| \check x\|^2 \le b \|x\| \cdot \| g(x)\|$ applies only to vectors $x\in \mathbb R^d_+$ with $\|\check x\|=o(\|x\|)$ for $\| x\| \to \infty$. Since also $\check x=0$ for $x \in \overline r$,  the condition defines a certain vicinity of the ray $\overline r$ (depending on $g$). Outside this region the relation between $g$ and $\sigma^2$ stays arbitrary.

For our second result on divergence of $(X_n)_{n\ge 0}$ we first rule out an evident case. We assume 
\begin{align}
 \exists u>0 \ : \ \mathbf P\big(\, X_n \to  X_\infty   \text{ with }u \le \|X_\infty\| < \infty\big)=0   \ . \tag{A2}\label{A2} 
\end{align}
Moreover we strengthen \eqref{smallo2} to the assumption
\begin{align} 
\exists \kappa >1/\delta \ : \  \sigma(x) = O(\|x\|\log^{-\kappa} \|x\|) \text{ for } \|x\|\to \infty \   , \tag{A3} 
\label{A3}
\end{align} 
where $\delta$ is as in assumption \eqref{A1}. 
\paragraph{Theorem 2.} {\em Let \eqref{A1} to \eqref{A3} be fulfilled and let $\varepsilon >0$. Assume that  for every $b>0$ there exists some $a>0$ such that for $x \in \mathbb R_+^d$
\begin{align} \| x\| \ge a \ , \ \| \check x\|  \le b \sigma (x) \quad \Rightarrow \quad \ell x \cdot \ell g(x) \ge\frac{1+ \varepsilon}2 \sigma^2(x) \ . 
\label{sigmacond}
\end{align}
Then there is a real number $v \ge0$ such that
\[ \mathbf P\big( \, \limsup_n \|X_n\| \le v \text{ or } \|X_n \| \to \infty \big) = 1\ . \]
If also for every $c>0$ there is a $n\in \mathbb N_0$ such that $\mathbf P(\| X_n \| > c ) > 0 $, then
\[ \mathbf P(\|X_n\| \to \infty ) > 0  \quad \text{ and }\quad \mathbf P\Big(\frac{X_n}{ \| X_n \| }\to \frac r{\|r\|}\ \Big|\ \|X_n\|\to \infty \Big)=1 \ . \]}

\mbox{}\\
Again we recover for $d=1$ the corresponding result from \cite{ke1}. Due to \eqref{A3}  it is now the condition $\| \check x\| \le b \sigma(x)$ giving the vicinity of the ray $\overline r$, where $g(x)$ and $\sigma^2(x)$ are interrelated. 

\paragraph{Remarks.} Let us comment on the assumptions of Theorem 2.

1. Obviously \eqref{A2} is also a necessary requirement in Theorem 2. Typically it is easily checked  in concrete examples. For Markov chains with a countable discrete state space $S \subset \mathbb R^d_+$ it says that away from zero there are no absorbing states. In the general case there is the following criterion: \eqref{A2} holds if $\ell g(x)$ is uniformly bounded away from zero on sets of the form $\{ x \in \mathbb R^d_+: u \le \ell x \le u+1 \}$ with $u >0$ sufficiently large. For the proof of this claim adopt the arguments at the end of section 2  in \cite{ke1} to the process $(\ell X_n)_{n \ge 0}$.

2.  Assumption \eqref{A3} cannot be weakened substantially in our general context. This follows from example C, Section 3 in \cite{ke1}. We note that \eqref{A3} is weaker than the corresponding assumption in \cite{ke1} for the 1-dimensional case. 

3. Remarkably, condition \eqref{sigmacond} cannot be relaxed in our general context. It is not enough to require \eqref{sigmacond} just for {\em some} $b>0$ as we shall see at the end of this paper by means of a counterexample. It is tempting to conjecture that condition \eqref{driftcond} cannot be weakened, too.
\qed

\mbox{}\\
So far we have not specified any choice of the norm $\|\ \|$ on $\mathbb R^d$. This was not necessary so far, since as is well-known all norms on a finite dimensional Euclidean space are equivalent, and one easily convinces oneself that all our conditions or statements involving norms are preserved if one passes to an  equivalent norm. Thus, in examples one may work with the most convenient one, e.g. the $l_1$- or $l_2$-norm. For our proofs these norms are not appropriate. We shall utilize a norm specificially suited for our purposes. This norm is introduced in section~2. The proofs of the theorems are then presented in section 3 and 4. They use ideas from \cite{ke1} and \cite{la} and are based on the construction of Lyapunov functions of the form
\[   l_{\alpha, \beta, \gamma,j}(x)=(1+ \gamma x_j/\ell x)\frac{\| \check x\|^2}{(\ell x)^2} (\log \ell x)^{-\beta-1} + \alpha (\log \ell x)^{-\beta} \]
with $x=(x_1, \ldots,x_d)^T \in \mathbb R^d_+$, $1  \le j \le d$, $\alpha >0$, $\gamma \ge 0$ and either $\beta=-1$ or $\beta >0$.  Section 5 contains the counterexample.

For notational convenience we use the symbol $c$ for a positive constant which may change its value from line to line.

\section{A useful norm}
Let us briefly put together the facts on matrices which we are going to use. Recall that $M$ is a primitive matrix with Perron-Frobenius eigenvalue 1 and corresponding left and right eigenvectors $\ell$ and $r$. Then as is well-known from Perron-Frobenius theory (see \cite{se}) 
\[  \max \{|\eta| : \eta \text{ is an eigenvalue of } M- r\ell \} < 1 \ .  \]
This maximum is called the {\em spectral radius} of the matrix $M-r\ell$. It follows from matrix theory (see \cite{ho}, Lemma 5.6.10) that one can construct a matrix norm  \mbox{$|\hspace{-1.5pt}|\hspace{-1.5pt}|  \  |\hspace{-1.5pt}|\hspace{-1.5pt}|$} on the space of all $d\times d$ matrices such that
\[ \rho:= |\hspace{-1.5pt}|\hspace{-1.5pt}|  M-r\ell   |\hspace{-1.5pt}|\hspace{-1.5pt}|  <1 \ . \]
From this matrix norm we obtain (see \cite{ho}, Theorem 5.7.13) a functional $\| \ \|$ on $\mathbb R^d$ via 
\[  \| x \| := |\hspace{-1.5pt}|\hspace{-1.5pt}| C_x |\hspace{-1.5pt}|\hspace{-1.5pt}| \ , \quad x \in \mathbb R^d\ , \]
where $C_x$ denotes the $d\times d$ matrix having all columns equal to $x$. $\| \ \|$ is a norm, since the properties of norms transfer from \mbox{$|\hspace{-1.5pt}|\hspace{-1.5pt}|  \  |\hspace{-1.5pt}|\hspace{-1.5pt}|$}  directly to $\| \ \|$. This is the norm we are going to work with in the sequel. It has the property 
\begin{align} \| Ax\| \le |\hspace{-1.5pt}|\hspace{-1.5pt}| A |\hspace{-1.5pt}|\hspace{-1.5pt}| \cdot \|x\| 
\label{matrixnorm}
\end{align}
for $x \in \mathbb R^d$ and  any $d\times d$ matrix $A$. Indeed  $C_{Ax} = AC_x$ and the property $|\hspace{-1.5pt}|\hspace{-1.5pt}| C_{Ax} |\hspace{-1.5pt}|\hspace{-1.5pt}|\le|\hspace{-1.5pt}|\hspace{-1.5pt}| A |\hspace{-1.5pt}|\hspace{-1.5pt}|\cdot |\hspace{-1.5pt}|\hspace{-1.5pt}| C_x |\hspace{-1.5pt}|\hspace{-1.5pt}|$  of matrix norms gives the claim.
In particular 
\begin{align}
\| (M-r\ell) x\| \le \rho \| x\| \ .
\label{contraction}
\end{align}
Thus $M-r\ell$ induces a contraction in the norm $\| \ \|$. 

By equivalence of norms we may change from $\| \ \|$ to any other norm. In particular there is a constant $\lambda< \infty$ such that
\begin{align} \| \check x\| \le \lambda \, \ell x \quad \text{for all } x \in \mathbb R_+^d \ . 
\label{estimate1}
\end{align}
To see this observe that from the inequality \eqref{matrixnorm} we have  $\| \check x\| \le \gamma \| x\| $ with $\gamma =$ \mbox{$|\hspace{-1.5pt}|\hspace{-1.5pt}|  I-r\ell  |\hspace{-1.5pt}|\hspace{-1.5pt}|$}.   Also $\| x\|' := \ell_1|x_1| + \cdots + \ell_d|x_d|$ defines a norm on $\mathbb R^d$, since $\ell_i>0$ for all $i=1, \ldots,d$. Thus by equivalence of norms we arrive at \eqref{estimate1}.

In order to apply these results to our process $(X_n)_{n \ge0}$ note that we have $(I-r\ell)M=M- r\ell=M(I-r\ell)$ and $\ell \check X_n =0$, thus
\begin{align*} \check X_{n+1} &= (I-r\ell)(MX_n+ g(X_n)+ \xi_{n+1}) \\
&= (M-r\ell)\check X_n + (I-r\ell)(g(X_n)+ \xi_{n+1}) \ .
\end{align*}
From \eqref{contraction} to \eqref{estimate1} it follows that 
\begin{align} \| \check X_{n+1}\| \le \rho \| \check X_n\| + c\, \ell g(X_n)+ c\, \|\xi_{n+1}\| 
\label{estimate2} 
\end{align}
for some $c < \infty$. (Here we need that $g(x)$ has only non-negative components.) 
Further observe that for any $\mu >0$ and $a,b \ge 0$ we have
\begin{align} (a+b)^2 \le (1+ \mu)  a^2 + (1+ \mu^{-1})  b^2 \ .
\label{estimate0}
\end{align}
Applying this estimate twice to the right-hand side of \eqref{estimate2} we obtain for any $\mu >0$
\begin{align}  \| \check X_{n+1}\|^2 \le (1+\mu)\rho^2 \|\check X_n\|^2 + c\, (\ell g(X_n))^2+ c\, \|\xi_{n+1}\|^2 
\label{estimate3}
\end{align}
with a suitable $c< \infty$.

\section{Proof of Theorem 1}
First observe that if we replace $X_n$ by $\overline X_n:= X_n + r$ for all $n \ge 0$ then equations \eqref{gleichung} and \eqref{variance} as well as assumption (A1) still hold, if $g(x)$ and $\sigma^2(x)$ are  replaced by $\overline g(x):=g(x-r)$ and $\overline \sigma^2(x)=:\sigma^2(x-r)$. Note that the assumptions \eqref{smallo} and \eqref{smallo2}  are not affected if $g$ and $\sigma^2$ are substituted by $\overline g$ and $\overline \sigma^2$, and the same holds true for the conditions formulated in Theorem 1 if one replaces $\varepsilon$ by $\varepsilon/2$. Thus without loss of generality we may assume $\ell X_n \ge 1$ for all $n \ge 0$ throughout the proof. Then for any $\alpha >0$ 
\[ L_n := \frac{\|\check X_n\|^2}{(\ell X_n)^2} + \alpha \log \ell X_n \ , \quad n \in \mathbb N_0 \ ,  \]
is a sequence of non-negative random variables. We show that for large $\alpha$ it possesses a supermartingale property. The proof uses the following estimate, where $I(A)$ denotes the indicator variable of an event $A$.

\paragraph{Lemma 1.} {\em For all $t>0,$ $ h > -t$ and $\eta >0$ }
\[\log (t+h) \le \log t + \frac ht - \frac 1{2(1+\eta)} \frac {h^2}{t^2} I(h \le \eta t) \ .\]

\begin{proof}
See formula (2) in \cite{ke1}. 
\end{proof}

\paragraph{Lemma 2.}
{\em If $\alpha$ is chosen large enough, then there is a number $s>0$   such that
\begin{align*}
\ell X_n \ge s \quad \Rightarrow \quad \mathbf E[ L_{n+1} \mid \mathcal F_n] \le L_n \text{ a.s.}
\end{align*}}

\begin{proof}
Since $\ell M=\ell$ we have the equation
\begin{align} \ell X_{n+1} = \ell X_n +\ell g(X_n)+ \ell \xi_{n+1} \ .
\label{gleichung2}
\end{align}
Thus $\ell \xi_{n+1} \ge -\mu\,\ell X_n$ implies $ \ell X_{n+1} \ge (1-\mu)\, \ell X_n$. 
Together with  \eqref{estimate3} and \eqref{estimate1} this entails
\begin{align} \frac{\|\check X_{n+1}\|^2}{(\ell X_{n+1})^2} \le\ \frac{ (1+\mu) \rho^2 \|\check X_n\|^2 + c\, (\ell g(X_n))^2+ c\, \|\xi_{n+1}\|^2 }{ (1-\mu)^2(\ell X_n)^2} 
+ \lambda^2 I(\ell\xi_{n+1} < -\mu \, \ell X_n) \label{estimate7}
\end{align}
for some sufficiently large $c<\infty$. Now $\rho < 1$, thus, if   $\mu$ is sufficiently close to 0, 
\begin{align*} \frac{\|\check X_{n+1}\|^2}{(\ell X_{n+1})^2} \le (1-\mu) \frac{\|\check X_n\|^2}{(\ell X_n)^2}+  c\frac{\, (\ell g(X_n))^2+ \, \|\xi_{n+1}\|^2 }{(\ell X_n)^2}+ \lambda^2 \frac{(\ell \xi_{n+1})^2}{\mu^2 (\ell X_n)^2} \ .
\end{align*}
In view of  \eqref{A1}, if we further enlarge $c$,
\begin{align} 
\mathbf E\Big[ \frac{\|\check X_{n+1}\|^2}{(\ell X_{n+1})^2}\ \big| \ \mathcal F_n\Big] \le (1-\mu) \frac{\|\check X_n\|^2}{(\ell X_n)^2} + c \frac{(\ell g(X_n))^2+ \sigma^2(X_n)}{ (\ell X_n)^2} \text{ a.s.}
\label{estimate4}
\end{align}

Next from \eqref{gleichung2} and Lemma 1 (with $t=\ell X_n+\ell g(X_n)$ and $h=  \ell \xi_{n+1}$) for $\eta >0$
\begin{align*}
\log \ell X_{n+1}\le \log (\ell X_n + \ell g(X_n)) + \frac{\ell\xi_{n+1}}{\ell X_n} 
- \frac{(\ell \xi_{n+1})^2}{2(1+\eta)(\ell X_n)^2}I\big(\ell\xi_{n+1} \le \eta(\ell X_n + \ell g(X_n))\big) \ .
\end{align*}
Taking the concavity of the log-function into account we get
\begin{align*}
\log \ell X_{n+1} &\le \log \ell X_n + \frac{\ell g(X_n)}{\ell X_n}+ \frac{\ell\xi_{n+1}}{\ell X_n} 
- \frac{(\ell \xi_{n+1})^2}{2(1+\eta)(\ell X_n)^2}+ \frac{(\ell \xi_{n+1})^2}{(\ell X_n)^2}I( \ell\xi_{n+1} > \eta \, \ell X_n )\ .
\end{align*}
Using \eqref{variance}, \eqref{A1} and the Markov inequality and choosing $\eta$ sufficiently small it follows
\begin{align*}
\mathbf E[ \log \ell X_{n+1} \mid \mathcal F_n] \le \log \ell X_n + \frac{\ell g(X_n)}{\ell X_n} - \frac{(1-\varepsilon/3)\sigma^2(X_n)}{2(\ell X_n)^2} + c \frac{\sigma^p(X_n)}{(\ell X_n)^p} \text{ a.s.}
\end{align*}
with some $c<\infty$. Because of \eqref{smallo2}  there is a number $s>0$ such that for $\ell X_n \ge s$
\begin{align}
\mathbf E[ \log \ell X_{n+1} \mid \mathcal F_n] \le \log \ell X_n + \frac{\ell g(X_n)}{\ell X_n} - \frac{ (1-\varepsilon/2)\sigma^2(X_n)}{2(\ell X_n)^2}  \text{ a.s.}
\label{estimate5}
\end{align}

Now combining \eqref{estimate4} and \eqref{estimate5} and using \eqref{smallo}  we get
\begin{align*}
\mathbf E[L_{n+1} \mid \mathcal F_n] \le L_n - \mu \frac{\| \check X_n\|^2}{(\ell X_n)^2}+ (\alpha+c) \frac{\ell g(X_n)}{\ell X_n} 
- \Big( \frac {1-\varepsilon/2}2 \alpha-c\Big) \frac{\sigma^2(X_n)}{(\ell X_n)^2} \text{ a.s.}
\end{align*}
for $\ell X_n \ge s$ and $s$ sufficiently large. If we let $\alpha \ge 6c/\varepsilon-c$ we arrive at
\begin{align*}
\mathbf E[L_{n+1} \mid \mathcal F_n] \le L_n - \mu \frac{\| \check X_n\|^2}{(\ell X_n)^2}+ (\alpha+c) \Big(\frac{\ell g(X_n)}{\ell X_n}  -   \frac {1-\varepsilon}2 \frac{\sigma^2(X_n)}{(\ell X_n)^2}\Big) \text{ a.s.}
\end{align*}
for $\ell X_n \ge s$.
We are now ready for  the conclusion: 

If $(\alpha +c)\ell g(X_n)\cdot \ell X_n \le \mu \|\check X_n\|^2$ , then obviously $\mathbf E[L_{n+1} \mid \mathcal F_n] \le L_n$ a.s. for $\ell X_n \ge s$.

If on the other hand $ \mu \|\check X_n\|^2\le (\alpha +c)\,\ell g(X_n)\cdot \ell X_n $ then by equivalence of norms there is a $b<\infty$ such that $\|\check X_n\|^2 \le b\, \|g(X_n)\|\cdot \|X_n\|$. Now the assumption of Theorem 1 comes into play, and again $\mathbf E[L_{n+1} \mid \mathcal F_n] \le L_n$ a.s., if only $\ell X_n$ is large enough. Thus the claim of the lemma follows.
\end{proof}

\noindent
We complete  the proof of Theorem 1 now as in \cite{ke1}. Suppose that the event $\|X_n\| \to \infty$ has positive probability. Then the same holds for the event $L_n \to \infty$,  and there is natural number $N$ such that $\mathbf P(E)>0$ for the event
\[ E=\{ \inf_{n \ge N} L_n \ge s, L_n \to \infty \} \ . \]
Define the stopping time
\[ T_N:= \min\{ n \ge N : L_n < s \} \ . \]
In view of Lemma 2 the process $(L_{n \wedge T})_{n \ge N}$ is a supermartingale. It is non-negative and thus a.s. convergent. However, on the event $E$ we  have $T_N= \infty$ and $L_n \to \infty$ and consequently $L_{n \wedge T}\to \infty$. This contradicts the assumption $\mathbf P(E)>0$, and the proof is finished.

\section{Proof of Theorem 2}

Here we may replace $X_n$ by $X_n+3r$. Therefore without loss of generality we assume $\ell X_n \ge 3$ for all $n \in \mathbb N_0$. Now we consider the processes $L=L^{\alpha,\beta,\gamma,j}$ given by
\[ L_n=L_n^{\alpha,\beta,\gamma,j}:= \frac{(1+\gamma X_{n,j}/\ell X_n)\|\check X_n\|^2}{(\ell X_n)^2(\log \ell X_n)^{\beta+1}} + \alpha (\log \ell X_n)^{-\beta}  \ , \ n \in \mathbb N_0\ , \]
with the $j$th component $X_{n,j}$ of $X_n$, $1 \le j \le d$, and with $\alpha, \beta >0$ and $\gamma \ge 0$. For convenience we only treat the case $2< p \le 3$ in which the following estimate is valid (in the case $p >3$ further terms of the Taylor expansion have to be considered).

\paragraph{Lemma 3.} {\em Let $\beta >0$ and $2 < p \le 3$. Set $f(t):= (\log t)^{-\beta}$. Then there is a constant $c<\infty$ such that for all $t \ge 3$ and $h > 3-t$}
\[f(t+h) \le f(t)+ f'(t)h + \frac 12 f''(t)h^2 +  \frac{c|h|^p}{ (\log t )^{\beta+1} t^p} + I(h \le - t/2) \ . \]

\begin{proof}
See formula (6) in \cite{ke1}. 
\end{proof}

\paragraph{Lemma 4.} {\em Let $0<\beta< \kappa \delta -1$ and $\gamma \ge 0$ such that $(1+\gamma/\ell_j)\rho^2 < 1$. Then, if $\alpha$ is sufficiently large, there is a real number $s>0$  such that
\begin{align*}
\ell X_n \ge s \quad \Rightarrow \quad \mathbf E[ L_{n+1}^{\alpha,\beta,\gamma,j} \mid \mathcal F_n] + \frac{\sigma(X_n)^p}{(\ell X_n)^p}\le L_n^{\alpha,\beta,\gamma,j}  \text{ a.s.}
\end{align*}}

\begin{proof}
We proceed similarly as in the proof of Lemma 2. Here instead of \eqref{estimate7} we have  the estimate
\begin{align*} &(1+\gamma X_{n+1,j}/\ell X_{n+1})\frac{\|\check X_{n+1}\|^2}{(\ell X_{n+1})^2(\log \ell X_{n+1})^{\beta+1}} \\ \mbox{} &\mbox{}\qquad\qquad \le  (1+\gamma/\ell_j)(1+\gamma X_{n,j}/\ell X_n) \frac{ (1+\mu) \rho^2 \|\check X_n\|^2 + c\, (\ell g(X_n))^2+ c\, \|\xi_{n+1}\|^2 }{ (1-\mu)^2(\ell X_n)^2(\log \ell X_n+ \log (1-\mu))^{1+\beta}} 
 \\ &\qquad\qquad\qquad\mbox{} 
 + \lambda^2 I(\ell\xi_{n+1} < -\mu \, \ell X_n) 
\end{align*}
By assumption on $\gamma$ and for $\mu>0$ sufficiently small this implies
\begin{align}
&\mathbf E\Big[ \frac{(1+\gamma X_{n+1,j}/\ell X_{n+1})\|\check X_{n+1}\|^2}{(\ell X_{n+1})^2(\log \ell X_{n+1})^{\beta+1}}\ \big|\ \mathcal F_n\Big] \label{estimate6} \\ 
&\ \mbox{}\qquad \qquad\le (1-\mu) \frac{(1+\gamma X_{n,j}/\ell X_{n})\|\check X_n\|^2}{(\ell X_n)^2(\log \ell X_n)^{\beta+1}}+ c \frac{(\ell g(X_n))^2+ \sigma^2(X_n)}{ (\ell X_n)^2 (\log \ell X_n)^{\beta +1}} + c \frac{\sigma^p(X_n)}{(\ell X_n)^p} \text{ a.s.} \notag
\end{align}
with some $c< \infty$. 
 
 Next from  Lemma 3  with $t=\ell X_n$ and $h= \ell g(X_n)+ \ell \xi_{n+1}$, from \eqref{estimate0} and \eqref{gleichung2}  and from $\ell g(X_n) \ge 0$
\begin{align*}
f(\ell X_{n+1} )&\le f(\ell X_n)+ f'(\ell X_n)(\ell g(X_n)+ \ell \xi_{n+1} )\\
&\qquad \mbox{} +\frac 12 f''(\ell X_n)((1+ \mu) (\ell \xi_{n+1})^2 + (1+\mu^{-1})(\ell g(X_n))^2) \\
&\qquad \mbox{} + c \frac{ (\ell g(X_n))^p + |\ell \xi_{n+1}|^p}{(\log \ell X_n)^{\beta+1} (\ell X_n)^p} +I( \ell \xi_{n+1} \le - \ell X_n/2)
\end{align*}
for a suitable $c>0$. Since $f''(t) \sim \beta(\log t)^{-\beta - 1} t^{-2}$ for $t \to \infty$,
\begin{align*}
\mathbf E[ f(\ell X_{n+1} )\mid \mathcal F_n]\le f(\ell X_n) &- \beta \frac {\ell g(X_n)}{ (\log \ell X_n)^{\beta+1} \ell X_n}
+ \frac \beta 2 \frac{ (1+ 2\mu) \sigma^2(X_n) + c(\ell g(X_n))^2}{ (\log \ell X_n)^{\beta+1} (\ell X_n)^2}\\
& \mbox{} + c \frac{ (\ell g(X_n))^p + \sigma^p(X_n)}{(\log \ell X_n)^{\beta+1} (\ell X_n)^p}+ c\frac{\sigma^p(X_n)}{(\ell X_n)^p}  \text{ a.s.}
\end{align*}
for $\ell X_n$ sufficiently large. Combining this estimate with \eqref{estimate6} and choosing $\alpha$ large enough we obtain in view of \eqref{smallo} and \eqref{smallo2} 
\begin{align*}
\mathbf E[ L_{n+1}  \mid \mathcal F_n]+ \frac{\sigma^p(X_n)}{(\ell X_n)^p}  \le\ L_n &- \mu\frac{ \|\check X_n\|^2}{(\ell X_n)^2(\log \ell X_n)^{\beta+1}} + ((\alpha+1) c+1) \frac{ \sigma^p(X_n)}{ (\ell X_n)^p}  \\
& \mbox{} + \frac{\alpha\beta} {(\log \ell X_n)^{\beta+1}} \Big( \frac{1+3\mu}2 \frac{\sigma^2(X_n)}{(\ell X_n)^2}-  (1-\mu)\frac {\ell g(X_n)}{\ell X_n} \Big)  \text{ a.s.}
\end{align*}
for $\ell X_n$ sufficiently large.  From \eqref{A2} we have for $0<\beta<\kappa\delta-1$
\[\frac{\sigma^p(x)}{(\ell x)^p} = O\big(\frac {\sigma^2(x)}{(\ell x)^2 (\log x)^{\kappa\delta}}\big) =o\big(\frac {\sigma^2(x)}{(\ell x)^2 (\log x)^{\beta+1}}\big)\quad \text{for } \|x\|\to \infty\ .\] 
Therefore for $0<\mu<1$ sufficiently small
\begin{align*}
\mathbf E[ L_{n+1}  \mid \mathcal F_n] + \frac{\sigma^p(X_n)}{(\ell X_n)^p}  \le\ & L_n -\mu  \frac{\|\check X_n\|^2}{(\ell X_n)^2(\log \ell X_n)^{\beta+1}}    \\
& \mbox{} \qquad + \frac{\alpha\beta(1-\mu)} {(\log \ell X_n)^{\beta+1}} \Big( \frac{1+ \varepsilon}2 \frac{\sigma^2(X_n)}{(\ell X_n)^2}-  \frac {\ell g(X_n)}{\ell X_n} \Big)  \text{ a.s.}
\end{align*}
if $\ell X_n$ is large enough. 
We come to the conclusion:

If $\| \check X_n\| \ge b\sigma(X_n)$ with some sufficiently large $b$, then the last estimate implies the claim $\mathbf E[ L_{n+1}\mid \mathcal F_n] + \sigma^p(X_n)/(\ell X_n)^p \le L_n $. If on the other hand $\| \check X_n\| \le b\sigma(X_n)$, then the assumption of Theorem 2 applies and again the claim follows.
\end{proof}

\mbox{}\\
For the proof of Theorem 2 we again construct a supermartingale, this time from $L=L^{\alpha,\beta,\gamma,j}$.  Observe that for some $s>0$ and for $m,m'>0$ and $t>s$ fulfilling
\[ \alpha (\log s)^{-\beta} \ge m> m'\ge (1+\gamma/\ell_j)\lambda^2(\log t)^{-\beta-1} + \alpha  (\log t)^{-\beta} \]
with $\lambda>0$ from formula \eqref{estimate1} we have
\begin{align*} 
L_n \le m \quad &\Rightarrow \quad \ell X_n \ge s \ , \\
L_n  \ge m' \quad &\Rightarrow \quad \ell X_n \le t \ .
\end{align*}
If we choose $\alpha$, $\beta$, $\gamma$ and $s$ as demanded in Lemma 4, then $(m \wedge L_n)_{n \ge 0}$ becomes a non-negative supermartingal, which thus is a.s. convergent. Then up to a null-event there arise three possibilities. Either $L_n \to 0$, then $\ell X_n \to \infty$. Or $\liminf_n L_n \ge m$, then $\limsup_n \ell X_n \le t$. Or else $L_n$ has a limit $0<L_\infty < m$, then $s\le \liminf_n \ell X_n <\infty$.

In order to transfer these alternatives to the process $(X_n)_{n \ge 0}$  we choose different $\beta_1,\beta_2>0$ and a $\gamma>0$ fitting the assumptions of Lemma 4. We consider the processes
\[ L^0:= L^{\alpha,\beta_1,0,1}\ ,\ L^1:= L^{\alpha,\beta_1, \gamma,1}\ ,\ \ldots\ ,\ L^d:= L^{\alpha,\beta_1, \gamma,d} \ ,\  L^{d+1}:=L^{\alpha, \beta_2,0,1}  \]
and for some $s,t,m>0$ the events
\begin{align*} &E:= \{ \ell X_n \to \infty\} \ , \ E' :=\ \{ \limsup_n \ell X_n \le t\}\ , \\ & E'' :=\ \bigcap_{i=0}^{d+1} \{s\le \liminf_n \ell X_n <\infty, L_n^i\to L_\infty^i \text{ with } 0<L_\infty^i < m \} 
\ .
\end{align*}
We let $\alpha,s,t$  large enough and $m$ small enough such that the above conclusion for $L=(L_n)_{n \ge 0}$ applies simultaneously to all processes $L^0, \ldots, L^{d+1}$. Then $\mathbf P(E\cup E'\cup E'')=1$. 

Let us show that $\mathbf P(E'')=0$ for $s$ sufficiently large. We have
\[ L^0_n= L^{d+1}_n (\log \ell X_n)^{\beta_2-\beta_1} \ .\]
Thus the sequence $\ell X_n$ is convergent on $E''$ with $s \le \lim_n \ell X_n < \infty$. This means that the random variables $\hat X_n =r\ell X_n$ converge on $E''$.
Next from the definition of $L^0$ it follows that the sequence $  \|\check X_n\| $ converges on the event $E''$ with some limit $Z$. If $Z=0$ then $\check X_n \to 0$, and we obtain that $X_n= \hat X_n+ \check X_n$ is convergent on $E''$.
If on the other hand $Z>0$, then we see from the convergence of $L^1_n, \ldots, L^d_n$ that the components $X_{n,1} , \ldots, X_{n,d}$ all converge on $E''$. Again we conclude that $X_n$ is a convergent sequence on the event $E''$. Let $X_\infty$ be the limit.

Now, given $u>0$, if we choose $s$ sufficiently large then from $s\le \lim_n \ell X_n<\infty$ on $E''$ we obtain    $u\le \|X_\infty\|< \infty$ by equivalence of norms. Therefore assumption \eqref{A2} may be applied and we obtain $\mathbf P(E'')=0$ and consequently $\mathbf P(E\cup E')=1$. By equivalence of norms this  translates into the first assertion of Theorem 2.

For the second assertion we switch back to the supermartingale $m\wedge L$ with $\gamma=0$. Let $c>t$ be such that 
\[\alpha (\log c)^{-\beta }+ \lambda^2 (\log c)^{-\beta-1} < \alpha (\log t)^{-\beta}\ . \] 
From the assumption of this assertion and by equivalence of norms there is a natural number $N$ such that $\mathbf P(\ell X_N > c) >0$. It follows 
\[\mathbf E[m \wedge L_N;\ell X_N >c ] < \alpha (\log t)^{-\beta}\mathbf P( \ell X_N >c)\ . \]From the supermartingale property of $m \wedge L$ and Fatou's Lemma
\[ \mathbf E[\lim_n m\wedge L_n;\ell X_N >c ] < \alpha (\log t)^{-\beta} \mathbf P( \ell X_N >c) \ . \]
If now $\mathbf P(E')=1$, then $\lim_n m \wedge L_n \ge \alpha (\log t)^{-\beta}$ a.s. which contradicts the last inequality. Therefore it follows $\mathbf P(E)>0$. This gives the second assertion.

For the last assertion we first show that
\begin{align}
\|\xi_{n+1}\| = o(\| X_n\|)  \text{ a.s. on the event } \| X_n\| \to \infty \ . 
\label{smallo3}
\end{align}
Define
\[ L_n':= L_n + \sum_{k =0}^{n-1} \frac {\sigma^p(X_k)}{(\ell X_k)^p }\]
and for a natural number $N$
\[ T_N:= \min \{ n \ge N : \ell X_n < s \} \ . \]
If again $\alpha, \beta, \gamma$ and $s$ are chosen in accordance with Lemma 4 then
$(L_{n \wedge T_N}')_{n \ge 0}$ is a non-negative supermartingal and thus a.s. convergent. It follows 
\[ \sum_{k=0}^\infty\frac {\sigma^p(X_k)}{(\ell X_k)^p } < \infty \text{ a.s. on the event } T_N = \infty \ . \]
Now in view of the first assertion of this theorem $\{ T_N = \infty\} \uparrow \{ \ell X_n \to \infty\} $ for $N\to \infty$, if only $s$ is sufficiently large. Therefore
\[ \sum_{k=0}^\infty\frac {\sigma^p(X_k)}{(\ell X_k)^p } < \infty \text{ a.s. on the event } \|X_n \| \to \infty \ .   \]
Because of \eqref{A1} and the Markov inequality this entails for every $\eta >0$
\[  \sum_{k=0}^\infty \mathbf P( \|\xi_{k+1}\| > \eta\, \ell X_k \mid \mathcal F_k) < \infty \text{ a.s. on the event } \|X_n \| \to \infty \ ,   \]
and the martingale version of the Borel-Cantelli Lemma (see \cite{du}, Theorem 5.3.2) implies \eqref{smallo2}.

Now from \eqref{estimate2}, \eqref{smallo} and \eqref{smallo3} we obtain that
\[ \| \check X_{n+1} \| \le \rho \| \check X_n\| + Y_n \quad \text{ with } Y_n = o(\|X_n \|) \text{ a.s. on } \| X_n\| \to \infty \ . \]
By induction 
\[ \| \check X_{n+1}\| \le \|\check X_0\| + \sum_{k=0}^n \rho^{n-k} Y_k \ . \]
Since $\rho < 1$ it follows 
\[ \| \check X_n\| = o(\|X_n\|) \text{ a.s. on the event } \|X_n\| \to \infty \ . \]
On the other hand $\hat X_n/\|\hat X_n\|= r/\|r\|$. This yields the last claim of Theorem 2.

\section{A counterexample}
We discuss an example in dimension $d=2$, which can be easily lifted to higher dimensions. In this section we use the $l_1$-norm $\| x\| := |x_1|+|x_2|$ for $x=(x_1,x_2)^T$. Let
\[ M= \frac 12 \begin{pmatrix} 1&1\\1&1 \end{pmatrix} \ , \quad \ r=\begin{pmatrix} 1\\ 1 \end{pmatrix}\ , \quad \ell= \frac 12 \big(1,1\big)\ .\]
Let $\overline g(t)$, $\overline \sigma(t)$, $t\ge 0$, be two   functions such that $\overline \sigma$ is differentiable and for $t \ge 0$
\[ \lim_{t\to \infty} \sigma'(t)= 0 \quad \text{ and }\quad \forall t > 0 : 0<  \overline g(t)\le \overline \sigma(t) \le t/2 \ ,\  |\sigma'(t)| < \frac 12\ . \] 
(For definiteness make $(0,0)^T$ an absorbing state.) Define for $x \in \mathbb R^2_+$
\[  \sigma (x):= \overline \sigma (\ell x) \ , \  g(x):= \begin{cases} \overline g(\ell x)r & \text{if } \|\check x\| \le \sigma(x) \\ (0,0)^T & \text{else .} \end{cases}  \]
Let $\chi_n$, $\zeta_n$, $n \ge 1$, be independent, $\mathbb R^2$-valued random variables  with
\begin{align*} \mathbf P(\chi_n=(1,1)^T)= \mathbf P(\chi_n=-(1,1)^T)= \mathbf P(\zeta_n=(1,-1)^T)= \mathbf P(\zeta_n=(-1,1)^T)= \frac 12 \ . 
\end{align*}
Define the Markov chain $X=(X_n)_{n \ge 0}$ inductively by $X_0=r$,
\[ \xi_{n+1} := \sigma (X_n) \chi_{n+1}+ \sigma(X_n) \zeta_{n+1}I(\|\check X_n\|\le \sigma (X_n) )  \]
and \eqref{gleichung}. 
Note that $M$ is the orthogonal projection on the subspace spanned by $r$. This together with the condition $\sigma (x) = \overline\sigma(\ell x) \le \ell x/2$ guarantees that the process $X$ never exists from the quadrant $\mathbb R^2_+$.  The conditions assumptions \eqref{smallo}, \eqref{variance}, \eqref{smallo2} and \eqref{A1} are fulfilled, and the same is true for \eqref{A2} and \eqref{A3} under mild conditions on $\overline g$ and $\overline \sigma$. However, due to the definition of $g(x)$, the condition \eqref{sigmacond} will never be satisfied for $b >1$, no matter how $\overline g$ and $\overline \sigma$ are chosen. We shall see that indeed the conclusion of Theorem 2 fails, even though \eqref{sigmacond} can be  achieved for  $b\le 1$ (but not all $b$). The reason is  that the process $X$ again and again leaves the region defined by the inequality $\|\check x\| \le \sigma(x)$.  

To prove this claim notice that from our assumptions for $t>0$
\[ \overline\sigma(t+ \overline g(t) \pm \overline \sigma(t) ) < \overline \sigma(t) + \tfrac 12 (\overline g(t)+ \overline \sigma(t)) \le 2 \overline \sigma(t) \ .\]
If now $\check X_n=0$ then from the definitions
\[ \ell X_{n+1} = \ell X_n + \overline g(\ell X_n) + \overline \sigma (\ell X_n) \ell \chi_{n+1} \quad \text{and} \quad \|\check X_{n+1} \| =\overline \sigma (\ell X_n) \|\zeta_{n+1}\| = 2\sigma (X_n)  \ .\]
From the previous inequality it follows $\sigma(X_{n+1})<2\sigma(X_n)$.  Thus $\sigma(X_{n+1}) < \|\check X_{n+1}\|$ and consequently from our definitions $ \check X_{n+2}=0$.  

Therefore, since we started with $\check X_0=0$, we have $ \check X_{2n}=0$ and $\| \check X_{2n+1}\| >\sigma (X_{2n+1}) $ for all $n \in \mathbb N_0$.
 Then $ \hat X_{2n}$, $ n \ge 0$, or (what amounts to the same thing) $\overline X_n:=\ell X_{2n}$, $n \ge 0$, is a Markov chain. Inserting our definitions we get
\[ \overline X_{n+1} =  \overline X_n + \overline g(\overline X_n) + \overline \xi_{n+1} \quad \text{ with } \quad \overline \xi_{n+1}:= \overline \sigma(\overline X_n) \ell \chi_{2n+1} + \overline \sigma (\ell X_{2n+1})\ell \chi_{2n+2} \ . \]
Letting $\overline {\mathcal F}_n := \mathcal F_{2n}$
\[ \mathbf E[ \overline \xi_{n+1} \mid \overline{\mathcal F}_n]= 0 \ , \ \mathbf E[ \overline \xi_{n+1}^2 \mid \overline{\mathcal F}_n]= \tau^2 (\overline X_n)  \]
with
\begin{align*} \tau^2(t) &= \overline \sigma^2(t)+ \mathbf E[ \overline \sigma^2(\ell X_1) \mid \ell X_0= t, \check X_0=0] \\ &
= \overline \sigma^2(t) + \frac 12 \overline \sigma^2(t +\overline g(t)+\overline \sigma(t))+ \frac 12 \overline \sigma^2(t +\overline g(t)-\overline \sigma(t))
\end{align*}
From our assumptions
\[ \tau^2(t) \sim 2\overline \sigma^2(t) \quad \text{for } t\to \infty \ .\]
Thus we are ready to apply our theorems (with $d=1$) to the process $\overline X=(\overline X_n)_{n \ge 0}$ 
and see that we have recurrence if $t \overline g(t) \le (1-\varepsilon) \overline \sigma^2(t)$ for large $t$. Note the the factor $1/2$ dropped out on the right-hand side. Thus there are cases, where the statement is false that there is transience for $t \overline g(t) \ge (1+\varepsilon) \overline \sigma^2(t)/2$. This shows that the assertion of Theorem 2 cannot be applied to the process $X$.

\end{document}